\documentclass[11pt]{article}

\usepackage{amsmath}
\usepackage{amsfonts}
\usepackage{amssymb}
\usepackage{amsthm}
\usepackage{graphics}
\usepackage{amscd}
\usepackage{graphicx,epsfig}
\usepackage{color}
\usepackage[all]{xy}
\usepackage{url}

\DeclareMathOperator{\Div}{div}

\renewcommand{\epsilon}{\varepsilon}

\newcommand{\boD}{\mathcal{D}}

\newcommand{\boC}{\mathcal{C}}

\newcommand{\R}{\mathbb{R}}

\newcommand{\Z}{\mathbb{Z}}

\renewcommand{\S}{\mathbb{S}}
\newcommand{\N}{\mathbb{N}}

\newcommand{\eps}{\varepsilon}

\newcommand{\der}[2]{\dfrac{\partial #1}{\partial #2}}
\newcommand{\dd}{\mathrm{d}}

\newcommand{\Ome}{\Omega}

\newtheorem{thm}{Theorem}

\newtheorem{prop}[thm]{Proposition}
\newtheorem{lem}[thm]{Lemma}

\renewcommand{\phi}{\varphi}
\newcommand{\dis}{\displaystyle}

\newtheorem*{thm*}{Theorem}
\newtheorem{prob}{Problem}

\theoremstyle{remark}

\newtheorem*{rem*}{Remark}

\newcounter{remark}

\newcounter{case}

\newcounter{construction}

\newcounter{fact}

\title{Adding one handle to half-plane layers}
\author{Laurent Mazet}
\date{}

\begin{document}

\maketitle

\begin{abstract}
In this paper, we build properly embedded singly periodic minimal
surfaces which have infinite total curvature in the quotient by their
period. These surfaces are constructed by adding a handle to the
toroidal half-plane layers defined by H.~Karcher. The technics that we
use is to solve a Jenkins-Serrin problem over a strip domain and to
consider the conjugate minimal surface to the graph.
\end{abstract}

\section*{Introduction}
In a preceding paper \cite{MRT}, the author with M.~Rodriguez and
M.~Traizet has contructed a familly of properly embedded singly
periodic minimal surfaces with an infinite number of
Scherk-ends. These surfaces are built as conjugate surfaces to
Jenkins-Serrin graphs over unbounded convex polygonal domains which
are neither a strip nor a half-plane. 

In this paper, we study the case where this unbounded convex polygonal
domain is a strip. If we do exactly the same construction as
in \cite{MRT} for a strip, we will prove that the minimal surface that
we obtain is actually doubly periodic. In fact, this give a new
construction of well-known examples of properly embedded doubly
periodic minimal tori with parallel ends. These surfaces are called
\emph{toro\"\i dal half-plane layers} by H.~Karcher who has built them 
in \cite{Ka}. In the classification of properly embedded doubly
periodic minimal tori with parallel ends made by J.~P\'erez,
M.~Rodriguez and M.~Traizet in \cite{PRT}, these surfaces are denoted
by $M_{\theta,\alpha,\frac{\pi}{2}}$ (see also \cite{Ro}). 

The aim of this paper is to modify the above example in such a way that
it loses its second period. The idea is to add one handle to this
minimal surface. F.~S.~Wei \cite{We} has been the first one to add
handles to half-plane layers but he can only do it in a periodic way
thus the surface is still doubly periodic. In \cite{RTW}, W.~Rossman,
E.~C.~Thayer and M.~Wohlgemuth have also added handles to the toro\"\i
dal half-plane layers in a periodic way. In fact we will prove:  
\begin{thm}\label{main}
There exists a properly embedded singly periodic minimal surface in
$\R^3$ whose quotient in $\R^2\times\S^1$ has genus $1$, an
infinite number of parallel Scherk ends and two limit ends.
\end{thm}

We notice that, recently, L.~Hauswirth and F.~Pacard \cite{HP} have
constructed examples of properly embedded minimal surface in $\R^3$
with two limit ends asymptotic to half Riemann surfaces. Their
construction consists in gluing a Costa-Hoffmann-Meeks surface of
small genus between two half Riemann surfaces.

Let us give some explanations on our construction. Recall how Scherk's
singly periodic surface can be built. From H.~Jenkins and J.~Serrin
\cite{JS}, there exists over the unit square a solution to the minimal
surface equation which takes the value $+\infty$ on two opposite sides
of the square and $-\infty$ on the other sides. The graph of this
solution is a minimal surface bounded by four vertical straight-lines
over the vertices of the square. This graph is a fundamental piece for
the Scherk's doubly periodic surface. The conjugate surface to this
graph is a minimal surface bounded by four horizontal symmetry curves
lying in two horizontal planes at distance $1$ from each other. By
reflecting about one of the two symmetry planes, we get a fundamental
domain for Scherk's singly periodic surface which has period
$(0,0,2)$.

This construction was generalized by H. Karcher \cite{Ka} and others
authors. In this paper we replace the unit square by a strip. We see
the boundary of the strip as the union of infinitly many unitary
edges. So we want to find a solution to the minimal surface equation
which take on the boundary the value $\pm\infty$ in alternating the
sign on every edge. In the second section, we build such a solution and
we describe the surface that we obtain when we consider the conjugate
surface to the graph: as we said above this surface was already built
by H.~Karcher and is called half-plane layer.

In the last section, we prove that we can add a handle to a half-plane
layer. This time, we solve a Dirichlet problem for the minimal surface
equation on a punctured strip. The solution is in fact multivaluated
in such a way that the graph of this solution is bounded by a vertical
straight-line over the removed point. The conjugate surface $\Sigma$ to
this graph is \textit{a priori} periodic but we prove that the period 
vanishes. Besides it is bounded by horizontal symmetry curves which lie
in the planes $\{z=0\}$ and $\{z=1\}$. By reflecting about one of this
two symmetry planes we obtain a fundamental piece of the surface
announced in Theorem \ref{main}. The vertical straight-line over the
removed point correspond in $\Sigma$ to a closed horizontal curve in
the handle.  

In the following section, we recall some tools that we use for the
study of the Dirichlet problem associated to the minimal surface
equation. 

The author would like to thank Martin Traizet for many helpful
discussions and for having talked him about this problem that was
asked by F.~S.~Wei.

\section{Preliminaries}

In this paper, we build minimal surfaces as graphs of functions. So in
this section, we recall some facts about the minimal surface
equation. 

Let $u$ be a function on a domain $\Ome\subset\R^2$. The graph of $u$
is a minimal surface if $u$ satisfies the minimal surface equation: 
\begin{equation*}\label{mse}
\Div\left(\frac{\nabla u}{\sqrt{1+|\nabla u|^2}}\right)=0
\tag{$*$}
\end{equation*}

In their article \cite{JS}, H.~Jenkins and J.~Serrin study the
Dirichlet problem associated to this partial differential equation in
bounded domain. They give necessary and sufficient conditions to
ensure existence of a solution.

One tools which is introduced is the conjugate $1$-form:
$$
\dd\Psi_u=\frac{u_x}{\sqrt{1+|\nabla u|^2}}\dd
y-\frac{u_y}{\sqrt{1+|\nabla u|^2}}\dd x
$$
$\dd\Psi_u$ is closed since $u$ satisfies \eqref{mse}. Then the
function $\Psi_u$ is locally defined and is called the conjugate
function. $\Psi_u$ has a geometric meaning since it is the third
coordinate of the conjugate minimal surface to the graph of $u$
expressed in the $x,y$ coordinates. If $v=\Psi_u$, we have $|\nabla
v|<1$, then $v=\Psi_u$ is $1$-Lipschitz continuous. Besides the
conjugate function $v$ is a solution of the maximal surface equation: 
\begin{equation*}\label{Mse}
\Div\left(\frac{\nabla v}{\sqrt{1-|\nabla v|^2}}\right)=0
\tag{$**$}
\end{equation*}

Any solution $v$ to \eqref{Mse} satisfies $|\nabla v|<1$. Then such a
solution extends continuously to the boundary. When the function $v$
can be writen $v=\Psi_u$, the computation of $v$ along the boundary is
given by lemmas in \cite{JS}. Let $v$ be a solution to \eqref{Mse} and
let us define the following $1$-form:
$$
\dd \Phi_v=\frac{v_y}{\sqrt{1-|\nabla v|^2}}\dd
x-\frac{v_x}{\sqrt{1-|\nabla v|^2}}\dd y
$$
$\dd\Phi_v$ is closed because of \eqref{Mse}. Besides the function
$\Phi_v$ is a solution of \eqref{mse}. If $u$ is a solution of
\eqref{mse} and $v=\Psi_u$, we have $u=\Phi_v$.

In this paper, we need to study the convergence of sequences
of solutions of \eqref{mse} or \eqref{Mse}. Because of the
correspondence $u\longleftrightarrow \Psi_u$ and $v\longleftrightarrow
\Phi_v$ the study is the same for the two equations. Let us recall
some points of the study of the convergnce of a sequence $(u_n)$ of
solutions of \eqref{mse}. Here the convergence we consider is the
$C^\infty$ convergence on every compact set. One can find proofs of all 
the facts below in \cite{Ma2,Ma3}.

Let $p$ be a point in $\Ome$ such that $|\nabla u_n|(p)\rightarrow
+\infty$. We can assume that $\nabla u_n/|\nabla u_n|(p)\rightarrow
\nu$ where $\nu$ is an unitary vector in $\R^2$. Let $L'$ be the
straight-line which passes by $p$ and is normal to $\nu$,  we then
denote by $L$ the connected component of $L'\cap \Ome$ which contains
$p$. Then for every point $q$ in $L$ we have $|\nabla u_n|(q)\rightarrow
+\infty$ and $\nabla u_n/|\nabla u_n|(q)\rightarrow \nu$. $L$ is
called a divergence line of $(u_n)$. Let us fix the orientation of $L$
such that $\nu$ points to the right-hand side of $L$. Then if
$q_1,q_2\in L$, the convergence of the gradient implies:
\begin{equation}\label{div1}
\lim_{n\rightarrow +\infty}\int_{[q_1,q_2]}\dd\Psi_{u_n}=|q_1q_2|
\end{equation}
where $|q_1q_2|$ is the distance between $q_1$ and $q_2$.

To prove the convergence of a subsequence of $(u_n)$, it suffices to
prove that there is no divergence line. Indeed, if it is the case, the
sequence $(|\nabla u_n|)$ is bounded at every point; then some
estimates on the derivatives implies that we can ensure the
convergence of a subsequence of $(u_n-u_n(p))$ for a point $p$ in
$\Ome$. Since we shall study only Dirichlet problem with infinite
boundary data, we notice that the vertical translation by $u_n(p)$
does not matter. We shall use some results that were developed in
\cite{Ma2} to prove that there is no divergence line.

Let us study a sequence $(v_n)$ of solutions of \eqref{Mse}. Since
$|\nabla v_n|<1$, we consider points $p$ where $\nabla v_n\rightarrow
\nu$ with $\nu$ a unitary vector of $\R^2$. The divergence line $L$
becomes the straight-line which is spanned by $\nu$ and, if
$q_1,q_2\in L$, Equation \eqref{div1} becomes
$$
\lim_{n\rightarrow +\infty}\int_{[q_1,q_2]}\dd v_n=|q_1q_2|
$$

\section{The first familly: the toro\"\i dal 
half-plane layers}\label{1ere}
In this section, we build a familly of doubly periodic properly
embedded minimal tori with parallel ends. These surfaces are called
half-plane layers and were defined by H.~Karcher \cite{Ka} (see also
\cite{Ro}). Here these surfaces are constructed as the conjugate
surface to a Jenkins-Serrin graph over an unbounded domain. 

\subsection{The Dirichlet problem}

We begin in precising the domains we shall consider.

Let $P$ be a polygon in $\R^2$, we say that $P$ is \emph{a convex
  unitary $2k$-gon} if $P$ is convex, has $2k$ edges and each one
  has unitary length. Let $P$ be a convex unitary $2k$-gon, $P$ is
  said to be special if $k\ge3$ and there exist $v,w\in\S^1$ such that
  two edges of $P$ are equal to $\pm v$ and all other edges are equal
  to $\pm w$. For more explanations on this definition, we refer to
  \cite{PT} and \cite{MRT}. 

Let $\Ome$ be a bounded polygonal domain, we say that $\Ome$ is a
\emph{unitary $2k$-polygonal domain} if the polygon $P=\partial\Ome$
is a convex unitary $2k$-gon. Let $\Ome$ be such a polygonal
domain. We denote by $e_1,\cdots,e_{2k}$ the $2k$ edges of
$\partial\Ome$ (the edges are numbered with respect to the orientation
of $\partial \Ome$).

One can consider the following Dirichlet problem :
\begin{prob}\label{prob1}
To find a solution $u$ of \eqref{mse} in $\Ome$ such that $u$ takes
the value $+\infty$ on $e_{2p}$ and $-\infty$ on $e_{2p-1}$ for
$p\in\{ 1,\cdots,k\}$
\end{prob}
It is known that this Dirichlet problem has a solution if and only if
$\partial\Ome$ is not special (see \cite{PT} and \cite{MRT}). Besides
the solution is unique.  

We want to study this Dirichlet problem when the domain $\Ome$
becomes a strip. We notice that the case when $\Ome$ becomes an
unbounded convex polygonal domain which is not a strip is studied in
\cite{MRT}. Let  $\Ome$ be a strip of width $y_0>0$. We can normalize 
$\Ome$ to be $\R\times (0,y_0)$. For every $n\in\Z$, let $a_0(n)$
(resp. $a_1(n)$) denote the point $(2n,0)$ (resp. $(2n+1,0)$). Let
$x_0\in[-1,1]$, then for every $n\in\Z$ we denote by $b_0(n)$
(resp. $b_1(n)$) the point $(x_0+2n-1,y_0)$
(resp. $(x_0+2n,y_0)$). Hence the Dirichlet problem we will study is: 
\begin{prob}
To find a solution $u$ of \eqref{mse} in $\Ome$ such that $u$ takes
the value $+\infty$ on $(a_0(n),a_1(n))$ and $(b_1(n),b_0(n+1))$ and
the value $-\infty$ on $(a_1(n),a_0(n+1))$ and $(b_0(n),b_1(n))$ for
every $n\in\Z$.  
\end{prob}

This Dirichlet problem is parametrized by the parameter $(x_0,y_0)$
which fixes the domain $\Ome$ and the boundary value. When $(x_0,y_0)$
describes $[-1,1]\times\R_+^*$, it is clear that we get all the
Dirichlet problem in a strip that can be seen as the limit of Problem
\ref{prob1}. Besides, the Dirichlet problem 
parametrized by $(-1,y_0)$ and $(1,y_0)$ are the same. Thanks to
Collin's example \cite{Co}, it is known that we do not have  
uniqueness of solutions to this Dirichlet problem.  Then to get
uniqueness, we add an other boundary condition. This condition has an
important geometrical meaning that we will explain soon. Thus for
$(x_0,y_0)\in [-1,1]\times \R_+$, we denote by $\boD(x_0,y_0)$ the 
following Dirichlet problem:
\begin{prob}[$\boldsymbol{\boD(x_0,y_0)}$]
To find a solution $u$ of \eqref{mse} in $\Ome$ such that $u$ takes
the value $+\infty$ on $(a_0(n),a_1(n))$ and $(b_1(n),b_0(n+1))$ and
the value $-\infty$ on $(a_1(n),a_0(n+1))$ and $(b_0(n),b_1(n))$ for
every $n\in\Z$ and 
\begin{equation}\label{equa1}
\int_{a_0(0)}^{b_1(0)}\dd\Psi_u=1
\end{equation}
\end{prob}
Let $u$ be a solution of this Dirichlet problem. Since $\Ome$ is
simply connected, the conjugate function $\Psi_u$ is well defined; we
fix $\Psi_u$ by $\Psi_u(a_0(0))=0$. So because of the value of $u$
on the boundary, for every $n$, if $p\in[a_1(n-1),a_1(n)]$,
$\Psi_u(p)=|a_0(n)p|$ (see \cite{JS}) ($|pq|$ denotes the euclidean
distance between $p$ and $q$). The graph of $\Psi_u$ on
$\R\times\{0\}$ is then serrate. Because of \eqref{equa1},
$\Psi_u(b_1(0))=1$ and, with the value of $u$ on the boundary, for
every $n\in\Z$ and $p\in[b_1(n-1),b_1(n)]$, $\Psi_u(p)=|b_0(n)p|$.

Before solving the Dirichlet problem $\boD(x_0,y_0)$, we need a lemma.
\begin{lem}\label{join}
Let $(x,y)$ be in $[-1,1]\times\R_+$ such that $x^2+y^2>1$. There
exists a convex polygonal path $\gamma$ in $[-1,+\infty)\times [0,y]$ that
  joins $(0,0)$ to $(x,y)$ and is composed of an odd number of unitary
  length edges.
\end{lem}
\begin{proof}
First there are two simple cases : $x=\pm 1$. If $x=1$, we begin
$\gamma$ by one edge that joins $(0,0)$ to $(1,0)$ then we joins
$(1,0)$ to $(1,y)$ by passing by the point $p=(t,y/2)$; $t$ is chosen
such that the distance between $p$ and $(1,0)$ is an integer. When
$x=-1$, this is the same idea : we end $\gamma$ by an edge that joins
$(0,y)$ to $(-1,y)$ and complete $\gamma$ as an isosceles triangle
(see Figure \ref{fig1}).

\begin{figure}[h]
\begin{center}
\resizebox{0.8\linewidth}{!}{\input{figanse1.pstex_t}}
\caption{\label{fig1}}
\end{center}
\end{figure}

Now, we assume that $x^2<1$. Let $n\in\N$ be a large integer and
$\theta>0$ be a small angle. We try to define $\gamma$ as follow : the 
$n$ first unitary edges of $\gamma$ joins $(0,0)$ to $(n\cos\theta,
n\sin\theta)$, the $n$ last unitary edges joins $(x+n,y)$ to $(x,y)$,
then, if the distance between $(n\cos\theta, n\sin\theta)$ and
$(x+n,y)$ is one, we can complete $\gamma$ by one unitary edge and
$\gamma$ will have $2n+1$ edges. So the idea is that we can choose
$n$ and $\theta$ such that this construction works and $\gamma$ has all
the expected properties (see Figure \ref{fig1}).

We want to chose $n$ and $\theta$ such that :
$$
\left(x+n(1-\cos\theta)\right)^2+\left(y-n\sin\theta\right)^2=1
$$
Then, with $t=1/n$, this equation becomes 
$$
(x^2+y^2-1)t^2+2\left(x(1-\cos\theta)- y\sin\theta\right)t+
2(1-\cos\theta)=0 
$$
For small $\theta>0$, the discrimant is $4\theta^2(1-x^2)+o(\theta^2)>
0$. Then one solution is:
$$
t=\theta\left(\frac{y+\sqrt{1-x^2}}{x^2+y^2-1}\right)+o(\theta)
$$
Hence there is small $\theta$ such that $1/t$ is an integer and we can
complete $\gamma$. Let us see that $\gamma$ satisfies the desired
properties. First:
$$
n\sin\theta= \frac{\sin\theta}{t}= \frac{x^2+y^2-1}{y+\sqrt{1-x^2}}
+o(1) 
$$
Since $x^2-1<0$, $n\sin\theta<y$ then $\gamma$ is in
$[-1,+\infty)\times[0,y]$. Besides 
$$
\frac{y-n\sin\theta}{(x+n)-n\cos\theta}= \frac{1}{x} \left( y-
\frac{x^2+y^2-1}{y+\sqrt{1-x^2}}\right)+o(1) 
$$
So when $x> 0$,
$\dis\frac{y-n\sin\theta}{(x+n)-n\cos\theta}>\tan\theta$ for small
$\theta$; hence $\gamma$ is convex. When $x\le 0$, $\gamma$ is clearly
convex.  
\end{proof}

We then can solve the Dirichlet problem.
\begin{thm}\label{dirich}
Let $(x_0,y_0)$ be in $[-1,1]\times \R_+$. Then there exists a
solution to $\boD(x_0,y_0)$ if and only if $x_0^2+y_0^2>1$. Besides,
when $x_0^2+y_0^2>1$, the solution is unique up to an additive
constant. 
\end{thm}

\begin{proof}
First, since $|\dd\Psi_u|<1$, \eqref{equa1} proves that the condition
$x_0^2+y_0^2>1$ needs to be satisfied for having a solution.

We now assume that $x_0^2+y_0^2>1$, we will prove the existence of a 
solution. We shall build a sequence $(\Ome_n)$ of bounded convex
unitary $2k$-polygonal domain in $\Ome$ such that the parallelogram
$a_0(n)b_1(n) b_0(-n+1)a_1(-n)$ is included in $\Ome_n$. We build
$\Ome_n$ as follows. From Lemma \ref{join}, there exists a convex
polygonal arc $\gamma$ that joins $a_0(n)$ to $b_1(n)$ with an odd
number of unitary edges in $[2n-1,+\infty)\times[0,y_0]$. Let $s$ be
  the symmetry with respect to the middle point of $[a_1(0),b_1(0)]$
  (\textit{i.e.} the point $((x_0+1)/2,y_0/2)$). The polygonal arc
  $s(\gamma)$ joins $b_0(-n+1)$ to $a_1(-n)$ in $(-\infty,-2n+1]\times
[0,y_0]$. We then define $\Ome_n$ as the polygonal domain bounded by
the following polygonal arc : the segment $[a_1(-n),a_0(n)]$, the arc
$\gamma$, the segment $[b_1(n),b_0(-n+1)]$ and the arc $s(\gamma)$. By
construction, $\Ome_n$ is a convex polygonal domain with an
even number of unitary edges and the parallelogram $a_0(n)b_1(n)
b_0(-n+1)a_1(-n)$ is included in $\Ome_n$. Besides we notice that
$\Ome_n$ is symmetric with respect to the middle point of
$[a_1(0),b_1(0)]$. 

Since $x_0^2+y_0^2>1$, the polygonal domain $\Ome_n$ is not special
(see \cite{MRT})
so there exists a solution $u_n$ of \eqref{mse} in $\Ome_n$ that
takes the value $+\infty$ on $[a_0(0),a_1(0)]$ and alternately the
values $+\infty$ and $-\infty$ on each edge of $\partial\Ome_n$. Let
$\Psi_n$ be the conjugate function to $u_n$, we fix $\Psi_n$ by
$\Psi_n(a_0(0))=0$. By construction, the number of edges between
$a_0(0)$ and $b_1(0)$ in $\partial\Ome_n$ is odd then
$\Psi_n(b_1(0))=1$. Besides, by maximum principle, $0\le\Psi_n\le 1$ in
$\Ome_n$. We restrict the function $u_n$ to the parallelogram 
$a_0(n)b_1(n)b_0(-n+1)a_1(-n)$ and we study the convergence of the
sequence $(u_n)$ on this increasing sequence of domain. Let us notice
that on each edge of the parallelogram which is included in
$\partial\Ome$, the function $u_n$ takes the value that we prescribe
in our Dirichlet problem $\boD(x_0,y_0)$.

Let us study the lines of divergence of this sequence on the limit
domain $\Ome$. Because of the value on the boundary the end points of
a line of divergence must be vertices of $\Ome$ (see Lemma A.3 in
\cite{Ma2}). Let $L$ be a line of divergence and let $p$ and $q$ in
$L$. Since $L$ is a divergence line, we have for a subsequence: 
\begin{equation}\label{div}
\lim_{n\rightarrow+\infty}\left|\Psi_n(p)-\Psi_n(q) \right|=|pq|
\end{equation}
Since for every $n$, $0\le\Psi_n\le 1$, this implies that $|pq|\le
1$. Then $L$ must have two end-points: one in $\{y=0\}$ and one in 
$\{y=y_0\}$. On each vertex of $\Ome$, $\Psi_n$ takes the value $0$ or
$1$. If we apply \eqref{div} with $p$ and $q$ the end-points of
$L$, we can conclude that one end-point is $a_i(k)$ and the other one
is $b_{1-i}(l)$ with $i\in\{0,1\}$ and $1=|a_i(k)b_{1-i}(l)|$. But
since $x_0^2+y_0^2>1$, the distance between one $a_i(k)$ and one
$b_{1-i}(l)$ is greater than $1$. Then the sequence has no divergence
line.

Then a sub-sequence of $(u_n)$ converges in $\Ome$ to a solution of
\eqref{mse}. Since the sequence $u_n$ takes the good boundary value
on each edges for every $n$, the function $u$ takes the expected
boundary value. Besides:
$$
\int_{[a_0(0),b_1(0)]}\dd\Psi_u=\lim_{n\rightarrow+\infty}
\int_{[a_0(0),b_1(0)]}\dd\Psi_n=\Psi_n(b_1(0))-\Psi_n(a_0(0)) =1 
$$
This end the proof of the existence.

We finish by the uniqueness proof. As we explain above, if $u$ and
$u'$ are two solutions to $\boD(x_0,y_0)$, the boundary value of
$\Psi_u$ and $\Psi_{u'}$ are the same. Then $\Psi_u$ and $\Psi_{u'}$
are two bounded solutions of the maximal surface equation with the
same boundary value ($\Psi_u$ and $\Psi_u'$ are bounded since every
point in $\Ome$ is at a bounded distance from the boundary).  The
uniqueness result in \cite{Ma1} (Corollary 4) implies that $u-u'$ is
constant. 
\end{proof}

Let us make one remark. In the construction, we have $0\le\Psi_n\le 1$,
then the solution $u$ satisfies $0\le\Psi_u\le 1$.

\subsection{The half-plane layers}

Using the solution to $\boD(x_0,y_0)$, we are then able to build a doubly
periodic properly embedded minimal tori with parallel ends.

Let $(x_0,y_0)$ be in $[-1,1]\times \R_+$ such that
$x_0^2+y_0^2>1$. Let $u$ be the solution to $\boD(x_0,y_0)$ given by
Theorem \ref{dirich}. First, because of the boundary value of $u$, the
boundary of the graph of $u$ is composed of vertical straight lines
over each vertex of $\Ome$.

Let $t$ by the translation by the vector $(-2,0)$. The function $u\circ
t$ which is defined in $\Ome$ is also a solution to
$\boD(x_0,y_0)$. Hence there exists a constant $k\in\R$ such that
$u\circ t=u+k$. The graph of $u$ is then periodic with $(2,0,k)$ as
period. Besides $u\circ t=u+k$ implies $\Psi_u\circ t=\Psi_u+k'$;
the boundary values show that $k'=0$.

Let $\Sigma$ denote the conjugate minimal surface to the graph of
$u$. $\Sigma$ is the conjugate surface to a graph over a convex
domain, by R.~Krust Theorem, $\Sigma$ is also a graph; it is then
embedded. Since  $0\le\Psi_u\le 1$, $\Sigma$ is included in
$\R^2\times[0,1]$. 

Let us study the boundary of $\Sigma$. In the neighborhood of each
vertex $p$ of $\Ome$, the graph of $u$ is bounded by a vertical
straight-line. In $\Sigma$, the conjugate of this vertical
straight-line is an horizontal symmetry curve. This curve is in the
$\{z=0\}$ plane when $p=a_0(k)$ or $p=b_0(k)$ (because $\Psi_u(p)=0$)
and it is in the $\{z=1\}$ plane when $p=a_1(k)$ or $p=b_1(k)$. The
boundary of $\Sigma$ is then composed of all these horizontal curves,
they are drawn in Figure \ref{fig2}.

\begin{figure}[h]
\begin{center}
\resizebox{0.8\linewidth}{!}{\input{figanse2.pstex_t}}
\caption{\label{fig2}}
\end{center}
\end{figure}

Since the graph of $u$ is periodic, the minimal surface $\Sigma$ is
also periodic. Let $X_1^*$, $X_2^*$ and $X_3^*$ be the three functions
on $\Ome$ that give the three coordinates of $\Sigma$. We have
$X_3^*=\Psi_u$. The period of $\Sigma$ is given by:
$$
\left(X_1^*,X_2^*,X_3^*\right)(2,y_0/2)-
\left(X_1^*,X_2^*,X_3^*\right)(0,y_0/2) 
$$
Since $\Psi_u\circ t=\Psi_u$, the third coodinate of the period is
zero. As $\Sigma$ is a graph, the horizontal part of the period does
not vanish: $\Sigma$ has a non-zero horizontal period.

Since $\Sigma$ is embedded in $\R^2\times[0,1]$ with boundary in
$\{z=0\}$ and $\{z=1\}$, we can extend $\Sigma$ into
$\widetilde{\Sigma}$ by symmetry with respect to all the horizontal
planes $\{z=n\}$ ($n\in\Z$). The surface $\widetilde{\Sigma}$ is then a
doubly periodic embedded minimal surface. The two periods are the
horizontal one that comes from $\Sigma$ and the vertical period
$(0,0,2)$ that comes from the horizontal symmetries.

\section{How to add one handle ?}
In this section, we shall construct a familly of singly-periodic
properly embedded minimal tori with an infinite number of parallel
ends and two limit ends. First let us consider one of the preceding
examples and consider it as an embedded minimal surface with infinite
total curvature in $\R^3/T$ with $T$ the vertical period (it has genus
$0$). Then the
idea to build our new familly is to had one handle to this periodic
minimal surface (see Figure \ref{fig3}). With this 
handle the surface loses its horizontal period so it has an infinite
number of ends. As we said in the introduction, F.~S.~Wei \cite{We},
W.~Rossman, E.~C.~Thayer and M.~Wohlgemuth \cite{RTW} have added in a
periodic way an infinite number of handles to the most symmetric
half-plane layers; we refer to their papers for pictures of such
surfaces. 

\begin{figure}[h]
\begin{center}
\resizebox{0.8\linewidth}{!}{\input{figanse3.pstex_t}}
\caption{\label{fig3}}
\end{center}
\end{figure}

As in the preceding section, to make this construction, we begin in
solving a Diriclet problem  and then we consider the conjugate minimal
surface to the graph.

\subsection{The Dirichlet problem}

Let $(x_0,y_0)$ be in $[-1,1]\times\R_+$, we denote by $\Ome$ the
polygonal domain associated to $(x_0,y_0)$ as in Section
\ref{1ere}. Let $c$ be the middle point of $[a_1(0),b_1(0)]$, the
coordinates of $c$ are $(\frac{x_0+1}{2},\frac{y_0}{2})$. 

We shall solve a Dirichlet problem for the maximal surface equation in
$\Ome\backslash\{c\}$. 

\begin{thm}\label{dirich2}
Let $(x_0,y_0)$ be in $[-1,1]\times\R_+$. We assume that
$(x_0+1)^2+y_0^2>4$. Then there exists a solution $v$ to \eqref{Mse}
with the following boundary values:
\begin{itemize}
\item if $p\in[a_1(n-1),a_1(n)]$, $v(p)=|a_0(n)p|$ 
\item if $p\in[b_1(n-1),b_1(n)]$, $v(p)=|b_0(n)p|$
\item $v(c)=1$
\end{itemize}
Besides the solution $v$ is unique.
\end{thm} 

\begin{proof}
First let us remark that the condition $(x_0+1)^2+y_0^2>4$ is
necessary since it says that the distance between $a_0(0)$ and $c$ is
greater than $1$ and $v(c)-v(a_0(0))=1$.

We notice that since a solution is bounded in the boundary and $\Ome$
is a strip a solution $v$ is bounded in $\Ome$; hence the uniqueness of
the solution is a consequence of Theorem 3 in \cite{Ma1}.  

Let us now prove the existence part of the theorem. Let $n$ be an
integer and consider the domain $\Ome_n$ that we have define in
Theorem \ref{dirich} proof. We also have the solution $u_n$ to
\eqref{mse} and its conjugate $\Psi_n$. Let us denote $\phi$ the
boundary value of $\Psi_n$ on $\partial\Ome_n$. 

Since $\Psi_n$ satisfies $|\nabla\Psi_n|<1$ in $\Ome_n$, we get
that for every $p,q\in\partial\Ome_n$: 
\begin{align}
|\phi(p)-\phi(q)|&\le d_{\Ome_n}(p,q) \text{ and } \label{cond1}\\
|\phi(p)-\phi(q)|&< d_{\Ome_n}(p,q) \text{ if }[p,q]\backslash
\partial\Ome_n \neq \emptyset \label{cond2}
\end{align}
Here, $d_{\Ome_n}$ denotes the intrinsic distance in $\Ome_n$, but
since $\Ome_n$ is convex this distance is the classical euclidean
distance in $\R^2$ (see also Theorem 1 in \cite{KM2}). 

Since $c$ is at a distance greater than $1$ from $a_0(0)$ and
$c\in[0,1]\times\R_+$, $c$ is at a distance greater than $1$ from
every point in $\partial\Ome_n$ where $\phi$ vanishes. This
implies that:
\begin{equation}\label{cond3}
\forall p\in\partial\Ome_n,\quad |\phi(p)-1|< d_{\Ome_n}(p,c)
\end{equation}

Then equations \eqref{cond1}, \eqref{cond2} and \eqref{cond3} implies
that there exists a solution $v_n$ of \eqref{Mse} in
$\Ome_n\backslash\{c\}$ which have $\phi$ as boundary value in
$\partial\Ome_n$ and $v_n(c)=1$; this result is a consequence of
Theorem 1 in \cite{KM2} and Theorem 4.1 in \cite{BS}. We remark that,
by maximum principle, we have $0\le v_n\le 1$. 

The solution $v$ will be constructed as the limit of the sequence
$(v_n)$. We consider the restriction of the function $v_n$ to the
parallelogram $a_0(n)b_1(n)b_0(-n+1)a_1(-n)$ minus the point $c$. This
increasing sequence of domains converges to $\Ome\backslash\{c\}$. We
notice that, as in Theorem \ref{dirich} proof, on each edges of
$\Ome$, $v_n$ takes the value which is prescribed in the Dirichlet
problem for $v$. 

Let us study the line of divergence of the sequence $(v_n)$. Let $L$
be such a line. Because of the boundary value of $v_n$ in
$\partial\Ome$, the end-points of $L$ can only be vertices of $\Ome$
or the point $c$. Let $p$ and $q$ be in $L$; since $L$ is a line of
divergence, we have for a subsequence:
\begin{equation}\label{equa2}
\lim_{n\rightarrow+\infty} |v_n(p)-v_n(q)|=|pq|
\end{equation}
Since $0\le v_n\le 1$, this equation implies that $|pq|\le 1$ then the
length of $L$ needs to be less than $1$: $L$ has two end-points. First
let us assume that one end-point of $L$ is $c$. If the other end-point
$p$ is one $a_1(k)$ or $b_1(k)$, $v_n(c)-v_n(p)=0$ then Equation
\eqref{equa2} implies $|cp|=0$  which is impossible. So the other
end-point $p$ is one $a_0(k)$ or $b_0(k)$. So $v_n(c)-v_n(q)=1$ and by
\eqref{equa2} $|cp|=1$; but $(x_0+1)^2+y_0^2>4$, then $|cp|>1$. Hence
one end-point of $L$ is in $\{y=0\}$ and the other is in
$\{y=y_0\}$. Since $x_0^2+y_0^2>1$, the same arguments as in Theorem
\ref{dirich} prove that such a line of divergence can not exist. Then
the sequence $(v_n)$ has no divergence line and for a subsequence
$(v_n)$ converges to a solution $v$ of \eqref{Mse} which takes the
expected boundary value on $\partial\Ome$ and $c$.
\end{proof}

We recall that $s$ denotes the symmetry with respect to the point
$c$. For every $n\in\N$ we have $s(a_0(n))=b_0(-n+1)$ and
$s(a_1(n))=b_1(-n)$. Then, if $v$ is the function given by Theorem
\ref{dirich2}, $v\circ s$ is also a solution to this Dirichlet
problem, then $v\circ s=v$. 

Let $t$ be the translation in $\R^2$ of vector $(-2,0)$. Let us study
the sequence $(v_n)=(v\circ t^n)$. Since $t(a_i(k))=a_i(k-1)$ and
$t(b_i(k))=b_i(k-1)$, the function $v_n$ is defined on
$\Ome\backslash\{t^{-n}(c)\}$ and its boundary value $v_n(p)$ is
$|a_0(k)p|$ in $[a_1(k-1),a_1(k)]$ and $|b_0(k)p|$ in
$[b_1(k-1),b_1(k)]$. The sequence of domains converges to $\Ome$. The
discussion in Theorem \ref{dirich} proof proves that $(v_n)$ has no
divergence line. So, for a subsequence, we can assume that $(v_n)$
converges to $v_\infty$ a solution in $\Ome$  with boundary value:
$$
v_\infty(p)=\begin{cases}
|a_0(k)p|\text{ in }[a_1(k-1),a_1(k)]\\
|b_0(k)p|\text{ in }[b_1(k-1),b_1(k)]
\end{cases}
$$
Let $u_0$ be the solution to $\boD(x_0,y_0)$ given by Theorem
\ref{dirich}, $v_\infty$ has the same boundary value as $\Psi_{u_0}$, then
by uniqueness of the solution $v_\infty=\Psi_{u_0}$. Besides this implies
that $\Psi_{u_0}$ is the only possible limit for a subsequence of $(v_n)$
then the sequence $(v_n)$, itself, converges to $\Psi_{u_0}$.

Let us make another remark. When the parameter $(x_0,y_0)$ moves in
its allowed domain, the domain $\Ome\backslash\{c\}$  and the boundary
value for $v$ change continuously. Hence the discussion about divergence
lines in Theorem \ref{dirich2} proof shows that the solution $v$
depends continuously of the parameter $(x_0,y_0)$.

\subsection{The minimal graph}

Let $(x_0,y_0)$ satisfies the hypotheses of Theorem \ref{dirich2} and
$v$ be the associated solution of \eqref{Mse} in $\Ome$. The $1$-form
$\dd\Phi_v$ locally defines a funtion $u$ which is a solution of
\eqref{mse}. In $\Ome\backslash\{c\}$, the function $u$ can be seen as
a multivaluated function \textit{i.e.} when we make a turn around $c$,
we add a constant $k$ to $u$. The graph of $u$ has a vertical period.

Because of the boundary value of $v$ along $\partial \Ome$, the
function $u$ takes  the value $+\infty$ on the edges $[a_0(n),a_1(n)]$
and $[b_1(n),b_0(n+1)]$ and the value $-\infty$ on $[a_1(n),a_0(n+1)]$
and $[b_0(n),b_1(n)]$ ($n\in\Z$) (see Lemma 4 in \cite{MRT}). A part of
the boundary of the graph of $v$ is then composed of vertical
straight-line over the vertices of $\Ome$. The last boundary part of
the graph is given by the behaviour near $c$. In fact the graph of $v$
is bounded by a a vertical straight-line over $c$: to see this, we
consider $u$ as a well defined function on the universal cover of a
neighborhood of $c$ and we apply Theorem 4.2 in \cite{Ma2}, $u$
satisfies the two hypotheses since the graph of $u$ has a vertical
period and, in $\Ome$, $v\le 1= v(c)$. 
Besides the equation $v\circ s=v$ implies that $u\circ s=u+k/2$; this
equation is written on the universal cover of $\Ome\backslash\{c\}$. 

The convergence $v\circ t^n\rightarrow \Psi_{u_0}$ implies that
$(u\circ t^n)$ converges to $u_0$. In a certain sense, it says that
the asymptotic behaviour of the graph of $u$ is given by the graph of
$u_0$. 

Besides, as $v$, the function $u$ depends continuously in the
parameter $(x_0,y_0)\in[-1,1]\times\R_+^*\cap\{(x_0+1)^2+y_0^2>4\}$. 

\subsection{The conjugate surface : half-plane layer with one handle}

In this subsection, we study the conjugate surface to the
graph of $u$. The surface that we obtain is the one that was announced 
in Theorem \ref{main}.  

Let $\Sigma$ be this conjugate. Since the graph of $u$ has a vertical
period, $\Sigma$ is also periodic; but in fact we have:
\begin{prop}
The periods of $\Sigma$ vanish.
\end{prop}

\begin{proof}
Let $\gamma$ be the circle of center $c$ and radius $\eps$ in $\Ome$;
$\gamma$ is a generator of the homotopy group of
$\Ome\backslash\{c\}$. Let $\dd X_1^*$, $\dd X_2^*$ and $\dd X_3^*$
the three coordinate $1$-forms in $\Ome\backslash\{c\}$ of the surface
$\Sigma$. These $1$-forms depend only in the first derivatives of $u$
which is why they are well defined on $\Ome\backslash\{c\}$. The period
of $\Sigma$ is then given by: 
$$
\int_\gamma\left(\dd X_1^*, \dd X_2^*, \dd X_3^*\right)
$$
Since $\dd X_3^*=\dd\Psi_u=\dd v$ the third coordinate of the period
vanishes. The equation $u\circ s =u+k/2$ implies that:
$$
s^*(\dd X_1^*,\dd X_2^*)=-(\dd X_1^*,\dd X_2^*)
$$
Since $s(\gamma)$ is in the same homotopy class of
$\pi_1(\Ome\backslash\{c\})$ as $\gamma$, we obtain.
\begin{align*}
\int_\gamma\left(\dd X_1^*, \dd X_2^*\right)=
\int_{s(\gamma)}\left(\dd X_1^*, \dd X_2^*\right)
&=\int_\gamma s^*\left(\dd X_1^*, \dd X_2^*\right)\\
&=\int_\gamma -\left(\dd X_1^*, \dd X_2^*\right) \\
&=-\left(\int_\gamma\left(\dd X_1^*, \dd X_2^*\right)\right)
\end{align*}
Then the two horizontal periods vanish. 
\end{proof}

$\Sigma$ has then the topology of an annulus. Besides since $u\circ
s=u+k/2$, the surface $\Sigma$ is symmetric with respect to a vertical
axis. We denote by $S$ this symmetry.

Since $0\le v\le 1$, the
surface $\Sigma$ is included in $\R^2\times[0,1]$. Let us study the
boundary of $\Sigma$. First, in a neighborhood of each vertex $p$ of
$\Ome$, the graph of $u$ is bounded by a vertical straight-line over
$p$. Then the conjugate of this line in $\Sigma$ is a horizontal
symmetry curve. Because of the value of $v$ on the vertices, if
$p=a_0(n)$ or $p=b_0(n)$, the curve is in the $z=0$ plane; if
$p=a_1(n)$ or $p=b_1(n)$, the curve is in the $z=1$ plane. The last
boundary component of $\Sigma$ comes from the conjugate of the
vertical line over $c$. This conjugate is a horizontal closed curve
$\Gamma$ in the $z=1$ plane: the curve is closed since the period
vanishes  and then it is the conjugate of a segment included in the
vertical line, $\Gamma$ is in $\{z=1\}$ since $v(c)=1$. Besides
$\Gamma$ is embedded since it is convex and has total curvature
$2\pi$. The convexity comes from the fact that it is the conjugate of
a vertical segment in the boundary of a graph and the total curvature
is $2\pi$ since this vertical segment correspond to one turn around
the point $c$ then the gauss map which is horizontal describes only
one time $\S^1$. Because of the symmetry of $\Sigma$, the Jordan curve
$\Gamma$ is symmetric with respect to $S$. All these boundary curves
are drawn in Figure \ref{fig3}. 

The last property of $\Sigma$ is the following.
\begin{prop}
The surface $\Sigma$ is embedded.
\end{prop}

\begin{proof}
Let $\Sigma_1$ be the conjugate of the part of the graph of $u$ which
is over $\Ome\cap\{y< y_0/2\}$. $\Sigma_1$ is a graph
over a domain $D_1$ since $\Ome\cap\{y< y_0/2\}$ is
convex then $\Sigma_1$ is embedded. Let $\Sigma_2$ be the conjugate of
the part of the graph of $u$ which is over
$\Ome\cap\{y> y_0/2\}$. $\Sigma_2$ is also a graph over
a domain that we denote by $D_2$. 

Let us study what happens over $y=y_0/2$. Let $\boC_-$ be the
projection in the horizontal plane $\{z=0\}=\R^2$ of the conjugate to
the curve in the graph of $u$ which is over
$(-\infty,(x_0+1)/2)\times\{y_0/2\}$, we parametrize $\boC_-$ by
$(-\infty,(x_0+1)/2)$. In the same way, let $\boC_+$ be the horizontal
projection of conjugate to the curve in the graph of $u$ which is over
$((x_0+1)/2,+\infty)\times\{y_0/2\}$, $\boC_+$ is parametrized by
$((x_0+1)/2,+\infty)$. 

Let us describe the boundary component of $D_1$ that corresponds to
the boundary component ${y=y_0/2}$ of $\Ome\cap\{y<
y_0/2\}$. This boundary component is composed of the union of
$\boC_-$, one half of the horizontal projection of $\Gamma$ and
$\boC_+$. The boundary component of $D_2$  associated to $\{y=y_0\}$
is the union of $\boC_-$, the other half of the projection of $\Gamma$ and
$\boC_+$. Since along $\boC_-$ and $\boC_+$ the domains $D_1$ and
$D_2$ are not on the same side of the boundary, the union of $\boC_-$,
$\boC_+$ and the horizontal projection of $\Gamma$ is embedded in
$\R^2$. 

Since the periods of $\Sigma$ vanish, there exist two functions
$X_1^*$ and $X_2^*$ on $\Ome\backslash\{c\}$ that give
the first two coordinates of $\Sigma$. We obtain:
$$
\boC_-'(x)= (\der{}{x} X_1^*,\der{}{x} X_2^*)(x,y_0/2)
$$
Let $n\in\N$, we have $\boC_-(x-2n)= (X_1^*,X_2^*) \circ
t^n(x,y_0/2)$. The convergence of $u\circ t^n$ to $u_0$ implies that
$(\der{}{x}X_1^*,\der{}{x}X_2^*)\circ t^n$ converges to
$(\der{}{x}{X_0}_1^*,\der{}{x}{X_0}_2^*)$ on $\Ome$ where ${X_0}_1^*$
and ${X_0}_2^*$ are the two first coordinates of the conjugate surface
to $u_0$. This implies that $||\boC_-'||$ is bounded in $\R_-$ then
$\boC_-$ is Lipschitz continuous. Besides:
\begin{align*}
\lim_{n\rightarrow+\infty}\boC_-(-2n+2)-\boC_-(-2n)
&=\lim_{n\rightarrow+\infty} \int_{-2n}^{-2n+2}\boC_-'(x)\dd x\\
&=\lim_{n\rightarrow+\infty} \int_0^2
(\der{X_1^*}{x},\der{X_2^*}{x})\circ t^n(x,y_0/2)\dd x\\
&=\int_0^2(\der{}{x}{X_0}_1^*,\der{}{x}{X_0}_2^*)(x,y_0/2)\dd x
\end{align*}
The last integral is the horizontal period of the conjugate surface to
the graph of $u_0$ then this vector is non zero. Since this limit does
not vanish and $\boC_-$ is Lipschitz continuous, we obtain 
$$\lim_{x\rightarrow-\infty}\boC_-(x)=\infty $$

Since the surface $\Sigma$ is symmetric, we have
$\boC_+(x)=S(\boC_-(x_0+1-x))$ then $\dis\lim_{x\rightarrow
+\infty}\boC_+(x)=\infty$. Let us compactify $\R^2$ into a sphere
$\S^2$, the curves $\boC_-$ and $\boC_+$ joins at infinity. Hence the
union of $\boC_-$, $\boC_+$ and $\Gamma$ cuts the sphere into three
connected components. Since $D_1$ and $D_2$ are connected and
$\boC_-$, $\boC_+$ and $\Gamma$ are in their boundary, each $D_i$ is
included in one of these three connected components. Over a
neighborhood of the point $\boC_-(-1)$, $\Sigma$ is graph then $D_1$
and $D_2$ do not lie in the same connected component of $\S^2
\backslash (\boC_-\cup\boC^+\cup\Gamma)$. Then $D_1$ and $D_2$ do not
intersect themselves. Hence $\Sigma_1$ and $\Sigma_2$ are disjointed and
so $\Sigma=\overline{\Sigma_1 \cup \Sigma_2}$ is embedded.
\end{proof}

Finally $\Sigma$ is an embedded minimal surface in $\R^2\times[0,1]$
with boundary in $\{z=0\}$ and $\{z=1\}$. Besides we can extend
$\Sigma$ into $\widetilde{\Sigma}$ by symmetry with respect to the
horizintal plane $\{z=n\}$. $\widetilde{\Sigma}$ is then an embedded
minimal surface in $\R^3$ with one vertical period $(0,0,2)$. The
quotient surface is a torus with an infinite numbers of parallel
Scherk ends. Theorem \ref{main} is then proved.

Let us also remark that because of the convergence of $(u\circ t^n)$
to $u_0$, we can say that $\Sigma$ has in a certain sense an
asymptotic behaviour given by the half-plane layer parametrized by
$(x_0,y_0)$. It implies, for example, that $\Sigma$ has infinite total
curvature. 

We notice that, since $v$ depends continuously on $(x_0,y_0)$, we can
say that, in a certain sense, $\Sigma$ depends continuously on
$(x_0,y_0)$. In particular, compact parts of $\Sigma$ depend
continuously.

\subsection{Two particular cases}
In this subsection, we study two particular range for the parameter
$(x_0,y_0)$. 

\subsubsection{$x_0=1$ and $y_0>0$}

In this case, the point $c$ is $(1,y_0/2)$. Besides the point $a_i(n)$
has the same abscissa as $b_i(n)$. Then, in the Dirichlet problem for
$v$, we have two symmetries for the boundary value $\phi$: 
\begin{enumerate}
\item $\phi(x,0)=\phi(x,y_0)$ for every $x\in\R$.
\item $\phi(x,y)=\phi(2-x,y)$ for every $x\in\R$ and $y\in\{y_0\}$.
\end{enumerate}
Because of the uniqueness in Theorem \ref{dirich2}, the solution $v$
has the same symmetries \textit{i.e.} $v(x,y)=v(x,y_0-y)$ and
$v(x,y)=v(2-x,y)$ for every $(x,y)\in\Ome\backslash\{c\}$. 

For the surface $\Sigma$, this implies that we have two vertical
orthogonal planes of symmetry. Let us observe the intersection of
$\Sigma$ with the plane that corresponds to the symmetry
$v(x,y)=v(2-x,y)$. This curve is the conjugate of the curve which is over
$\{x=1\}$ in the graph of $u$. This curve has two components which are
symmetric with respect to the other plane and each one joins a
boundary component in the plane $\{z=1\}$ to $\Gamma$ which is also in
$\{z=1\}$. Then, following the terminology of W.~Rossman, E.~C.~Thayer
and M.~Wohlgemuth in \cite{RTW}, the handle of the surface
$\widetilde{\Sigma}$ is a $'-'$ type handle.

We notice that combining the two planar symmetries we recover the
symmetry with respect to a vertical axis. Besides, to build the
surface $\Sigma$ we only need one fourth of this surface. This fourth
of a surface can be seen as the conjugate of the graph of $u$ over
$[1,+\infty) \times[0,y_0/2]$ which is convex. So this fourth of
$\Sigma$ is a graph over a domain, the symmetries implies that this
domain is included in an angular sector of angle $\pi/2$. Then, by
symmetry, the whole surface $\Sigma$ is a graph and is then
embedded. In this case, it is then easier to see the embeddedness. 

We remark that it seems possible (but difficult) to prove the
embeddedness of $\Sigma$ for all $(x_0,y_0)$ in using the embeddedness
in the particular case $x_0=1$ and the continuity in the parameter. 

\subsubsection{$x_0=-1$ and $y_0>2$}

In this case, the point $c$ is $(0,y_0/2)$ and the point $a_i(n)$ has
the same abscissa as $b_i(n+1)$. Then, we also have two symmetries for
the boundary value of $v$. Because of the uniqueness of the solution,
$v$ satisfies: $v(x,y)=v(x,y_0-y)$ and $v(x,y)=v(-x,y)$ for every
$(x,y)\in\Ome\backslash\{c\}$. 

So we have also two vertical orthognal planes of symmetry for the
minimal surface $\Sigma$. The intersection of $\Sigma$ with the plane
that corresponds to $v(x,y)=v(-x,y)$ is then the conjugate of the curve
which is over $\{x=0\}$ in the graph of $u$. This conjugate curve has
then two connected components which are symmetric with respect to the
other symmetry. Each component joins a boundary curve of $\Sigma$
which is in $\{z=0\}$ to $\Gamma\subset\{z=1\}$. Then, this times, the
handle of $\widetilde{\Sigma}$ is a $'+'$ type handle (see
\cite{RTW}).  

Thus, when $x_0$ goes from $1$ to $-1$, the surface
$\widetilde{\Sigma}$ is continuously deformed from a surface with a
$'-'$ type handle into a surface with a $'+'$ type handle. 

As in the preceding particular case, it is easy to prove the
embeddedness by using the symmetries.

\bigskip

\noindent Laurent Mazet

\noindent Laboratoire de Math\'ematiques et physique th\'eorique

\noindent Facult\'e des Sciences et Techniques, Universit\'e de Tours

\noindent Parc de Grandmont 37200 Tours, France.

\noindent E-mail: laurent.mazet@lmpt.univ-tours.fr

\end{document}